\documentclass[leqno]{llncs}

\usepackage{amsfonts}
\usepackage{amsmath}
\usepackage{amssymb}
\usepackage{comment}
\usepackage{expdlist}
\usepackage[dvips]{graphicx}
\usepackage{subfigure}
\usepackage{algorithmic}
\usepackage{float}
\usepackage{algorithm}

\newcommand{\set}[1]{\left \{ #1 \right \}}                     
\newcommand{\setst}[2]{\left \{ #1 \mid #2 \right \}}           
\newcommand{\abs}[1]{\left| #1 \right|}

\newtheorem{myclaim}{Claim}
\newtheorem{mydefinition}{Definition}

\renewcommand{\phi}{\varphi}
\renewcommand{\tilde}{\widetilde}

\renewcommand{\bar}{\overline}
\renewcommand{\epsilon}{\varepsilon}
\renewcommand{\vec}[1]{\overrightarrow{#1}}

\newcommand{\calF}{\mathcal{F}}
\newcommand{\calD}{\mathcal{D}}
\newcommand{\calT}{\mathcal{T}}

\newcommand{\nul}{\textsc{null}}
\newcommand{\true}{\textsc{true}}
\newcommand{\false}{\textsc{false}}
\newcommand{\symdiff}{\bigtriangleup}
\newcommand{\mydef}[1]{\mathop{\rm def}\nolimits ( #1 )}
\newcommand{\assign}{\Leftarrow}
\newcommand{\deltain}{\delta^{\rm in}}
\newcommand{\deltaout}{\delta^{\rm out}}

\newenvironment{numitem}{\refstepcounter{equation}\begin{enumerate}%
\item[(\theequation)]}{\end{enumerate}}

\newcommand{\refeq}[1]{(\ref{eq:#1})}                
\newcommand{\reffig}[1]{Fig.~\ref{fig:#1}}           
\newcommand{\reflm}[1]{Lemma~\ref{lm:#1}}            
\newcommand{\refsec}[1]{Section~\ref{sec:#1}}        
\newcommand{\refalgo}[1]{Algorithm~\ref{algo:#1}}    
\newcommand{\refclaim}[1]{Claim~\ref{cl:#1}}         

\renewenvironment{proof}{\par\noindent%
{\bf Proof.\par\nopagebreak}}{\unskip\nobreak\enskip$\square$\par\bigskip}
\newenvironment{proofof}[1]{\medskip\par\noindent%
{\bf Proof of #1.\par\nopagebreak}}{\unskip\nobreak\enskip$\square$\par\bigskip}

\makeindex

\begin{document}

\title{A Faster Algorithm for the \\ Maximum Even Factor Problem}

\institute
{
    Moscow State University
}

\author
{
   Maxim A. Babenko
   \thanks{
    Email: \texttt{max@adde.math.msu.su}.
        Supported by RFBR grant 09-01-00709-a.
   }
}

\maketitle

\begin{abstract}
    Given a digraph $G = (VG,AG)$, an \emph{even factor} $M \subseteq AG$ is a subset of arcs
    that decomposes into a collection of node-disjoint paths and even cycles.
    Even factors in digraphs were introduced by Geleen and Cunningham and generalize path
    matchings in undirected graphs.

    Finding an even factor of maximum cardinality in a general digraph
    is known to be NP-hard but for the class of \emph{odd-cycle symmetric} digraphs the problem
    is polynomially solvable.
    So far, the only combinatorial algorithm known for this task is due to Pap;
    it has the running time of $O(n^4)$ (hereinafter $n$ stands for the number of nodes in $G$).

    In this paper we present a novel \emph{sparse recovery} technique
    and devise an $O(n^3 \log n)$-time algorithm for finding a maximum cardinality
    even factor in an odd-cycle symmetric digraph.
\end{abstract}

\section{Introduction}
\label{sec:intro}

In \cite{CG-97} Cunningham and Geleen introduced the notion of \emph{independent path matchings}
and investigated their connection to the separation algorithms for the matchable set polytope,
which was previously studied by Balas and Pulleyblank~\cite{BP-89}.
Finding an independent path matching of maximum size was recognized as an intriguing example of a graph-theoretic
optimization problem that is difficult to tackle by the purely combinatorial means.
Two algorithms were given by Cunningham and Geleen: one relies on the ellipsoid method~\cite{CG-97},
and the other is based on deterministic evaluations of the Tutte matrix~\cite{CG-00}.
Later, a rather complicated combinatorial algorithm was proposed by Spille and Weismantel~\cite{SW-02}.

The notion of an \emph{even factor} was introduced as a further generalization
of path matchings in a manuscript of Cunningham and Geleen \cite{CG-01} (see also~\cite{cun-02}).
An even factor is a set of arcs that decomposes into a node-disjoint
collection of simple paths and simple cycles of even lengths.
Since cycles of length~2 are allowed, it is not difficult to see that finding a maximum matching
in an undirected graph~$G$ reduces to computing a maximum even factor in the digraph
obtained from~$G$ by replacing each edge with a pair of oppositely directed arcs.
On the other hand, no reduction from even factors to non-bipartite matchings is known.

Finding a maximum cardinality even factor is known to be NP-hard in general digraphs~\cite{cun-02}.
For the class of \emph{weakly symmetric digraphs} a min-max relation
was established by Cunningham and Geleen~\cite{CG-01} and by Pap and Szeg\H{o}~\cite{PS-04}.
Later it was noted by Pap \cite{pap-05} that these arguments hold for a slightly broader
class of \emph{odd-cycle symmetric digraphs}.
Takazawa and Kobayashi \cite{KT-09} pointed out a relation between even factors and jump systems and showed
that the requirement for a digraph to be odd-cycle symmetric is natural, in a sense.

The question of finding a combinatorial solution to the maximum even factor problem
in an odd-cycle symmetric digraph stood open for quite a while until Pap gave a direct $O(n^4)$-time algorithm~\cite{pap-05,pap-07}.
His method can be slightly sped up to $O(n^2(m + n \log n))$, as explained in~\refsec{simple_augment}.
(Hereinafter $n$ stands for the number of nodes in $G$ and $m$ denotes the number of arcs.)
To compare: the classical algorithm of Micali and Vazirani for finding a maximum non-bipartite matching,
which is a special case of the maximum even factor problem, runs in $O(mn^{1/2})$ time \cite{MV-80}.
It is tempting to design a faster algorithm for the maximum even factor problem
by applying the ideas developed for matchings (e.g. blocking augmentation~\cite{HK-73}).
There are, however, certain complications making even the bound of $O(mn)$ nontrivial.

To explain the nature of these difficulties let us briefly review Pap's approach
(a more detailed exposition will be given in~\refsec{simple_augment}).
It resembles Edmonds' non-bipartite matching algorithm and executes a series of iterations
each trying to increase the size of the current even factor~$M$ by one.
At each such iteration, a search for an augmenting path~$P$ is performed.
If no such path is found then~$M$ is maximum.
Otherwise, the algorithm tries to apply $P$ to $M$.
If no odd cycle appears after the augmentation then the iteration completes.
Otherwise, a certain contracting reduction is applied to~$G$ and $M$.

Hence, each iteration consists of \emph{phases} and the number of nodes in 
the current digraph decreases with each phase.
Totally there are $O(n)$ iterations and $O(n)$ phases during each iteration, which is quite
similar to the usual blossom-shrinking method. The difference is that during a phase
the reduction may completely change the alternating reachability structure so the next phase
is forced to start looking for $P$ from scratch.
(Compare this with Edmonds' algorithm were a blossom contraction changes the alternating forest
in a predicable and consistent way thus allowing this forest to be reused over the phases.)

In this paper we present a novel $O(n^3 \log n)$-time algorithm for solving the maximum even factor problem.
It is based on Pap's method but grows the alternating forest in a more careful fashion.
When a contraction is made in the current digraph the forest gets destroyed.
However, we are able to restore it by running a \emph{sparse recovery} procedure
that carries out a reachability search in a specially crafted digraph with $O(n)$ arcs
in $O(n \log n)$ time (where the $\log n$ factor comes from manipulations with data structures).

Our method seems applicable to a large variety of related problems.
In particular, the $O(mn^3)$-time algorithm of Takazawa \cite{T-08}
solves the weighted even factor problem in $O(mn^3)$ time and also involves
recomputing the alternating forest from scratch on each phase.

Similar effects of reachability failure are known to occur in the maximum $C_4$-free 2-factor problem \cite{pap-07}.
Also, adding matroidal structures leads to the \emph{maximum independent even factor problem}, which is solvable by the methods
similar to the discussed above, see~\cite{IT-07}.
All these problems can benefit from the sparse recovery technique.
Due to the lack of space we omit the details on these extensions.

\section{Preliminaries}
\label{sec:prelim}

We employ some standard graph-theoretic notation throughout the paper.
For an undirected graph~$G$, we denote its sets of nodes and edges by $VG$
and $EG$, respectively. For a directed graph, we speak of arcs rather than
edges and denote the arc set of $G$ by $AG$.
A similar notation is used for paths, trees, and etc.
We allow parallel edges and arcs but not loops.
As long as this leads to no confusion, an arc from $u$ to $v$ is denoted by $(u,v)$.

A path or a cycle is called \emph{even} (respectively \emph{odd}) if is consists
of an even (respectively \emph{odd}) number of arcs or edges.
For a digraph~$G$ a \emph{path-cycle matching} is a subset of arcs $M$
that is a union of node-disjoint simple paths and cycles in $G$. When $M$ contains
no odd cycle it is called an \emph{even factor}.
The \emph{size} of $M$ is its cardinality and the \emph{maximum even factor problem}
prompts for constructing an even factor of maximum size.

An arc $(u,v)$ in a digraph~$G$ is called \emph{symmetric} if $(v,u)$ is also present in~$G$.
Following the terminology from~\cite{pap-05}, we call~$G$
\emph{odd-cycle symmetric} (respectively \emph{weakly symmetric})
if for each odd  (respectively any) cycle $C$ all the arcs of~$C$ are symmetric.
As already noted in \refsec{intro}, the maximum even factor problem is NP-hard in general
but is tractable for odd-cycle symmetric digraphs.

Maximizing the size of an even factor $M$ in a digraph~$G$ is equivalent
to minimizing its \emph{deficiency} $\mydef{G,M} := \abs{VG} - \abs{M}$.
The minimum deficiency of an even factor in~$G$ is called
the \emph{deficiency} of~$G$ and is denoted by $\mydef{G}$.

For a digraph~$G$ and $U \subseteq VG$, the set of arcs entering (respectively leaving)~$U$
is denoted by $\deltain_G(U)$ and $\deltaout_G(U)$.
Also, we write $\gamma_G(U)$ to denote the set of arcs with both endpoints in~$U$ and
$G[U]$ to denote the subgraph of $G$ induced by $U$, i.e. $G[U] = (U, \gamma_G(U))$.
When, the digraph is clear from the context it is omitted from notation.

To \emph{contract} a set $U \subseteq VG$ in a digraph~$G$ means to replace nodes in $U$ by a single \emph{complex node}.
The arcs in $\gamma(VG - U)$ are not affected, arcs in $\gamma(U)$ are dropped, and
the arcs in $\deltain(U)$ (respectively $\deltaout(U)$) are redirected so as to enter (respectively leave) the complex node.
The resulting graph is denoted by $G / U$.
We identify the arcs in $G / U$ with their pre-images in $G$.
Note that $G / U$ may contain multiple parallel arcs but not loops.
If $G'$ is obtained from $G$ by an arbitrary series of contractions then
$G' = G / U_1 / \ldots / U_k$ for a certain family of disjoint subsets $U_1, \ldots, U_k \subseteq VG$
(called the \emph{maximum contracted sets}).

\section{Pap's Algorithm}
\label{sec:simple_augment}

Consider an odd-cycle symmetric digraph $G$.
The algorithm for finding a maximum even factor in $G$ follows the standard scheme of cardinality augmentation.
Namely, we initially start with the empty even factor~$M$ and execute a series of \emph{iterations}
each aiming to increase $\abs{M}$ by one.
Iterations call \textsc{Simple-Augment} routine that, given an odd-cycle symmetric digraph~$G$ and an even factor $M$ in $G$
either returns a larger even factor $M^+$ or $\nul$ indicating that the maximum size is reached.

\subsection{Augmentations}

Let us temporarily allow odd cycles and focus on path-cycle matchings in~$G$.
The latter are easily characterized as follows.
Construct two disjoint copies of $VG$: $V^1 := \setst{v^1}{v \in VG}$ and $V^2 := \setst{v^2}{v \in VG}$.
For each arc $a = (u,v) \in AG$ add the edge $\set{u_1, v_2}$ (\emph{corresponding} to~$a$).
Denote the resulting undirected bipartite graph by $\tilde G$.

Clearly, a path-cycle matching $M$ in $G$ is characterized by the following
properties: for each node $v \in VG$, $M$ has at most one arc entering $v$ and
also at most one arc leaving~$v$. Translating this to $\tilde G$ one readily sees
that $M$ generates a matching~$\tilde M$ in~$\tilde G$. Moreover, this correspondence
between matchings in $\tilde G$ and path-cycle matchings in $G$ is one-to-one.
A node $u^1$ (respectively $u^2$) in $\tilde G$ not covered by $\tilde M$ is called
a \emph{source} (a \emph{sink}, respectively).

\begin{algorithm}[t!]
    \caption{$\textsc{Simple-Augment}(G,M)$}
    \label{algo:simple_augment}
    \begin{algorithmic}[1]
        \STATE Search for an augmenting path~$P$ in $\vec G(M)$
        \IF{$P$ does not exist}
            \RETURN $\nul$
        \ELSIF{$P$ exists and is feasible}
            \RETURN $M \symdiff A(P)$
        \ELSE[$P$ exists but is not feasible]
            \STATE Put $M_i \assign M \symdiff A(P_i)$ for $i = 0, \ldots, k + 1$
            \STATE Find an index~$i$ such that $M_i$ is an even factor while $M_{i+1}$ is not
            \STATE Find the unique odd cycle $C$ in $M_{i+1}$
            \STATE $G' \assign G / C$, $M' \assign M_i / C$
            \STATE $\bar M' \assign \textsc{Simple-Augment}(G',M')$
            \IF[{$M'$ is maximum in $G'$}]{$\bar M' = \nul$}
                \RETURN $\nul$
            \ELSE[$M'$ is augmented in $G'$ to a larger even factor~$\bar M'$]
                \STATE Undo the contractions and transform $\bar M'$ to an even factor $M^+$ in $G$
                \RETURN $M^+$
            \ENDIF
        \ENDIF
    \end{algorithmic}
\end{algorithm}

Given a digraph~$G$ and a path-cycle matching~$M$ in~$G$ we turn $\tilde G$ into a digraph~$\vec G(M)$
by directing the edges $\set{u^1,v^2}$ corresponding to arcs $(u,v) \in M$ from $v^2$ to~$u^1$
and the other edges from $u^1$ to $v^2$.
\begin{mydefinition}
    A simple path in $\vec G(M)$ starting in a source node is called \emph{alternating}.
    An alternating path ending in a sink node is called \emph{augmenting}.
\end{mydefinition}
For an alternating path $P$ let $A(P)$ denote the set of arcs in $G$ corresponding to
the arcs of~$P$ in $\vec G(M)$.
Hereinafter $A \symdiff B$ denotes the symmetric difference of sets $A$ and $B$.
The next statements are well-known.

\begin{myclaim}
\label{cl:no_aug_path}
    If $\vec G(M)$ admits no augmenting path then $M$ is a path-cycle matching of maximum size.
\end{myclaim}

\begin{myclaim}
\label{cl:aug_path}
    If $P$ is an augmenting (respectively an even alternating) path in $\vec G(M)$ then $M' := M \bigtriangleup A(P)$ is path-cycle matching
    obeying $\abs{M'} = \abs{M} + 1$ (respectively $\abs{M'} = \abs{M}$).
\end{myclaim}

\begin{figure}[t!]
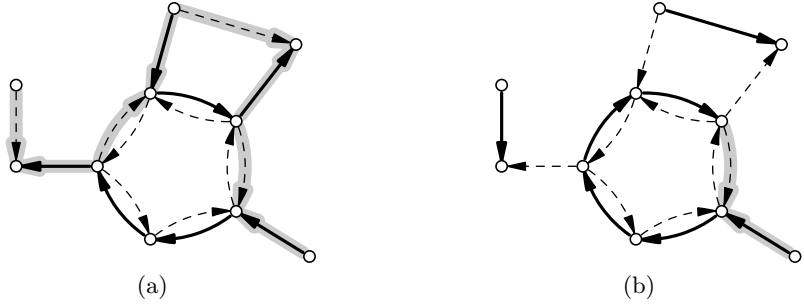

    \centering
    \subfigure[]{
        \includegraphics{pics/contraction.1}%
    }
    \hspace{2cm}
    \subfigure[]{
        \includegraphics{pics/contraction.2}%
    }
    \caption{
        Preparing for a contraction.
        Subfigure~(a): the arcs of~$M$ are bold and the grayed arcs correspond to path~$P_{i+1}$.
        Subfigure~(b): path~$P_i$ is applied, the arcs of~$M_i$ are bold, and the grayed arcs indicate the remaining part of~$P_{i+1}$.
    }
    \label{fig:contraction}
\end{figure}

The augmentation procedure (see~\refalgo{simple_augment}) constructs
$\vec G(M)$ and searches for an augmenting path~$P$ there. In case no such path exists,
the current even factor~$M$ is maximum by \refclaim{no_aug_path} (even in the broader class of path-cycle matchings),
hence the algorithm terminates.
Next, consider the case when $\vec G(M)$ contains an augmenting path~$P$.
\refclaim{aug_path} indicates how a larger path-cycle matching~$M'$ can be formed from~$M$,
however $M'$ may contain an odd cycle. The next definition focuses on this issue.
\begin{mydefinition}
    Let $P$ be an augmenting or an even alternating path.
    Then $P$ is called \emph{feasible}   if $M' := M \symdiff A(P)$ is again an even factor.
\end{mydefinition}

If~$P$ is feasible then \textsc{Simple-Augment} exits with the updated even factor $M \symdiff A(P)$.
Consider the contrary, i.e. $P$ is not feasible. Clearly, $P$ is odd, say it consists of $2k + 1$ arcs.
Construct a series of even alternating paths $P_0, \ldots, P_k$ where
$P_i$ is formed by taking the first $2i$ arcs of $P$ ($0 \le i \le k$). Also, put $P_{k+1} := P$
and $M_i := M \symdiff A(P_i)$ ($0 \le i \le k + 1$).

Then there exists an index $i$ ($0 \le i \le k$) such that $P_i$ is feasible
while $P_{i+1}$ is not feasible. In other words, $M_i$ is an even factor obeying $\mydef{G,M_i} = \mydef{G,M}$
and $M_{i+1}$ contains an odd cycle. Since~$M_i$ and $M_{i+1}$ differ by at most two arcs,
it can be easily showed that an odd cycle in $M_{i+1}$, call it $C$, is unique (see~\cite{pap-05}).
Moreover, $M_i$ \emph{fits} $C$, that is, $\abs{M_i \cap AC} = \abs{VC} - 1$ and $\deltaout(VC) \cap M_i = \emptyset$.
See \reffig{contraction} for an example.

It turns out that when an even factor fits an odd cycle
then a certain optimality-preserving contraction is possible.
As long as no confusion is possible, for a digraph $H$ and a cycle $K$ we abbreviate $H / VK$ to $H/K$.
Also, for $X \subseteq AH$ we write $X/K$ to denote $X \setminus \gamma_H(VK)$.
\begin{myclaim}[Pap~\cite{pap-05}]
\label{cl:contraction_opt}
    Let $K$ be an odd cycle in $H$ and $N$ be an even factor that fits~$K$.
    Put $H' := H / K$ and $N' := N / K$.
    Then $H'$ is an odd-cycle symmetric digraph and $N'$ is an even factor in~$H'$.
    Moreover, if $N'$ is maximum in~$H'$ then $N$ is maximum in~$H$.
\end{myclaim}

Note that $M$ is maximum in $G$ if and only if $M_i$ is maximum in $G$.
\textsc{Simple-Augment} contracts~$C$ in $G$.
Let $G' := G / C$ and $M' := M_i / C$ denote the resulting digraph and the even factor.
To check if $M'$ is maximum in $G'$ a recursive call $\textsc{Simple-Augment}(G',M')$ is made.
If $\nul$ is returned then $M'$ is a maximum even factor in $G'$, which implies by \refclaim{contraction_opt}
that the initial even factor $M$ was maximum in $G$. In this case \textsc{Simple-Augment} terminates returning $\nul$.

\begin{figure}[t!]
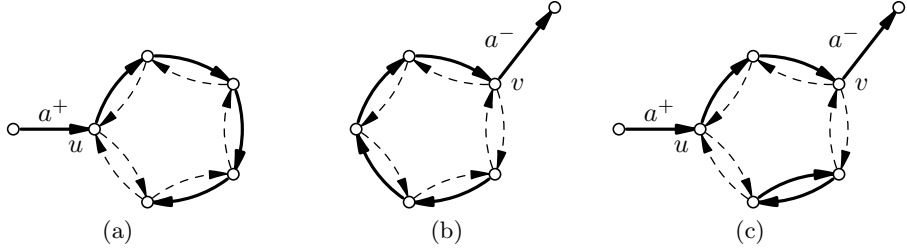

    \centering
    \subfigure[]{
        \includegraphics{pics/uncontraction.1}%
    }
    \hspace{1cm}
    \subfigure[]{
        \includegraphics{pics/uncontraction.2}%
    }
    \hspace{0.3cm}
    \subfigure[]{
        \includegraphics{pics/uncontraction.3}%
    }
    \caption{
        Some cases appearing in \reflm{uncontraction}. Digraph~$G$ and even factor~$N$ (bold arcs) are depicted.
    }
    \label{fig:uncontraction}
\end{figure}

Suppose that the recursive call was able to augment $M'$ to a larger even factor $\bar M'$ in~$G'$.
The algorithm transforms $\bar M'$ to an even factor in~$G$ with the help of the following statement from~\cite{pap-05}
(see~\reffig{uncontraction} for examples):
\begin{lemma}
\label{lm:uncontraction}
    Let $K$ be an odd cycle in an odd-cycle symmetric digraph~$H$.
    Put $H' := H / K$ and let $N'$ be an even factor in $H'$.
    Then there exists an even factor~$N$ in $H$ obeying $\mydef{H,N} = \mydef{H',N'}$.
\end{lemma}

\subsection{Complexity}

There are $O(n)$ iterations each consisting of $O(n)$ phases.
To bound the complexity of a single phase note that it takes $O(m)$ time to find an augmenting path~$P$
(if the latter exists). We may construct all the path-cycle matchings $M_0, \ldots, M_{k+1}$ and decompose
each of them into node-disjoint paths and cycles in $O(n^2)$ time. Hence, finding the index~$i$ and the cycle $C$
takes $O(n^2)$ time. Contracting~$C$ in~$G$ takes $O(m)$ time.
(We have already spent $O(m)$ time looking for $P$, so it is feasible to spend another $O(m)$ time
to construct the new digraph $G' = G/C$ explicitly.)
An obvious bookkeeping allows to undo all the contractions performed during the iteration
and transform the final even factor in the contracted digraph an even factor the initial digraph in $O(m)$ time.
Totally, the algorithm runs in $O(n^4)$ time.

The above bound can be slightly improved as follows.
Note that the algorithm needs an arbitrary index~$i$ such that $M_i$ is an even factor and $M_{i+1}$ is not,
i.e. $i$ is not required to be minimum.
Hence, we may carry out a binary search over the range $[0,k+1]$.
At each step we keep a pair of indices $(l,r)$ such that $M_l$ is an even factor while $M_r$ is not.
Replacing the current segment $[l,r]$ by a twice smaller one
takes $O(n)$ time and requires constructing and testing a single path-cycle matching $M_t$ where $t := \lfloor (l + r) / 2 \rfloor$.
This way, the $O(n^2)$ term reduces to $O(n \log n)$ and
the total running time of becomes $O(n^2 (m + n \log n))$.
The ultimate goal of this paper is to get rid of the $O(m)$ term.

\section{A Faster Algorithm}
\label{sec:fast_augment}

\subsection{Augmentations}

The bottleneck of \textsc{Simple-Augment} are the augmenting path computations.
To obtain an improvement we need better understanding of how these paths are calculated.
Similarly to the usual path-finding algorithms we maintain a directed out-forest $\calF$ rooted at the source nodes.
The nodes belonging to this forest are called \emph{$\calF$-reachable}.
At each step a new arc $(u,v)$ leaving an $\calF$-reachable node~$u$ is scanned and
either gets added to $\calF$ (thus making $v$ $\calF$-reachable) or skipped because $v$ is already $\calF$-reachable.
This process continues until a sink node is reached or no unscanned arcs remain in the digraph.

\begin{mydefinition}
    Let $G$ be a digraph and $M$ be an even factor in $G$.
    An \emph{alternating forest} $\calF$ for $M$ is a directed
    out-forest in $\vec G(M)$ such that:
    (i) the roots of $\calF$ are all the source nodes in $\vec G(M)$;
    (ii) every path from a root of~$\calF$ to a leaf of $\calF$ is even.
\end{mydefinition}
The intuition behind the suggested improvement is to grow~$\calF$ carefully
and avoid exploring infeasible alternating paths.
\begin{mydefinition}
    An alternating forest~$\calF$ is called \emph{feasible} if every even alternating or augmenting path in $\calF$ is feasible.
    An alternating forest~$\calF$ is called \emph{complete} if it contains no sink node and
    for each arc $(u,v)$ in $\vec G(M)$ if $u$ is $\calF$-reachable then so is $v$.
\end{mydefinition}

We replace \textsc{Simple-Augment} by a more sophisticated recursive \textsc{Fast-Augment} procedure.
It takes an odd-cycle symmetric digraph~$G$, an even factor $M$ in $G$, and an additional flag named $\textit{sparsify}$.
The procedure returns a digraph~$\bar G$ obtained from~$G$ by a number of contractions and an even factor~$\bar M$ in~$\bar G$.
Additionally, it may return an alternating forest~$\bar \calF$ for $\bar M$ in $\bar G$.
Exactly one of the following two cases applies:
\begin{numitem}
\label{eq:augment_breakthrough}
	$\mydef{\bar G,\bar M} = \mydef{G,M} - 1$ and $\bar\calF$ is undefined;
\end{numitem}
\begin{numitem}
\label{eq:augment_maximum}
	$\mydef{\bar G,\bar M} = \mydef{G,M}$,
        $\bar M$ is maximum in~$\bar G$,
        $M$ is maximum in~$G$, and
	$\bar\calF$~is a proper complete alternating forest for~$\bar M$ in~$\bar G$.
\end{numitem}

Assuming the above properties are true, let us explain how \textsc{Fast-Augment} can be used to perform a single augmenting iteration.
Given a current even factor~$M$ in~$G$ the algorithm calls $\textsc{Fast-Augment}(G,M,\true)$ and examines the result.
If $\mydef{\bar G,\bar M} = \mydef{G,M}$ then by \refeq{augment_maximum} $M$
is a maximum even factor in~$G$, the algorithm stops.
(Note that the forest~$\bar \calF$, which is also returned by \textsc{Fast-Augment}, is not used here.
This forest is needed due to the recursive nature of~\textsc{Fast-Augment}.)
Otherwise $\mydef{\bar G,\bar M} = \mydef{G,M} - 1$ by \refeq{augment_breakthrough};
this case will be referred to as a \emph{breakthrough}.
Applying \reflm{uncontraction}, $\bar M$ is transformed to an even factor $M^+$ in~$G$
such that $\mydef{G,M^+} = \mydef{\bar G,\bar M} = \mydef{G,M} - 1$.
This completes the current iteration.

\begin{algorithm}[t!]
    \caption{$\textsc{Fast-Augment}(G,M,\textit{sparsify})$}
    \label{algo:fast_augment}
    \begin{algorithmic}[1]
        \STATE Initialize forest~$\calF$ \label{line:forest_init}
        \FORALL{unscanned arcs $a = (u, v)$ such that $u^1 \in V\calF$} \label{line:loop_first}
            \STATE Mark $a$ as scanned \label{line:scan}
            \IF{$a \in M$ or $v^2 \in V\calF$}
                \STATE \textbf{continue for} \COMMENT{to line~2}
            \ENDIF
            \STATE $a_1 \assign (u^1,v^2)$
            \STATE Let $P_0$ be the even alternating path to $u^1$ in $\calF$
            \COMMENT{$P_0$ is feasible}
            \IF[single step]{$v^2$ is a sink}
                \STATE $P_1 \assign P_0 \circ a_1$  \label{line:single_step_first}
                \COMMENT{$P_1$ is augmenting}
                \IF{$P_1$ is feasible} \label{line:feasibility_check_1}
                    \RETURN $(G,M \symdiff A(P_1),\nul)$
                \ENDIF  \label{line:single_step_last}
            \ELSE[double step]
                \STATE Let $a_2 = (v^2, w^1)$ be the unique arc leaving $v^2$ \label{line:double_step_first}
                \COMMENT{$w^1 \notin V\calF$}
                \STATE $P_1 \assign P_0 \circ a_1 \circ a_2$
                \COMMENT{$P_1$ is even alternating}
                \IF{$P_1$ is feasible} \label{line:feasibility_check_2}
                    \STATE Add nodes $v^2$ and $w^1$ and arcs $a_1$, $a_2$ to $\calF$
                    \STATE \textbf{continue for} \COMMENT{to line~2}
                \ENDIF \label{line:double_step_last}
            \ENDIF

            \STATE{$M_0 \assign M \symdiff A(P_0)$, $M_1 \assign M \symdiff A(P_1)$} \label{line:contraction_start}
            \STATE{Find a unique odd cycle $C$ in $M_1$}
            \STATE{$G' \assign G / C$, $M' \assign M_0 / C$}

            \IF{$\textit{sparsify} = \false$}
                \RETURN $\textsc{Fast-Augment}(G',M',\false)$
            \ENDIF

            \STATE Construct the digraph $H'$
            \STATE $(\bar H', \bar M', \bar \calF) \assign \textsc{Fast-Augment}(H',M',\false)$ \label{line:recurse}
            \STATE Compare $V\bar H'$ and $VG$: let $Z_1, \ldots, Z_k$ be the maximum contracted sets
            and $z_1, \ldots, z_k$ be the corresponding complex nodes in~$\bar H'$
            \STATE $\bar G' \assign G / Z_1 / \ldots / Z_k$

            \IF{$\mydef{\bar G',\bar M'} < \mydef{G',M'}$}
                \RETURN $(\bar G',\bar M',\nul)$
            \ENDIF

            \STATE Unscan the arcs in $\bar G'$ that belong to~$M$ and the arcs that enter $z_1, \ldots, z_k$ \label{line:unscan}
            \STATE $G \assign \bar G'$, $M \assign \bar M'$, $\calF \assign \bar \calF$ \label{line:recovered}
        \ENDFOR \label{line:loop_last}
        \RETURN $(G,M,\calF)$
    \end{algorithmic}
\end{algorithm}

\medskip

Clearly, the algorithm constructs a maximum even factor correctly provided that
\textsc{Fast-Augment} obeys the contract. Let us focus on the latter procedure.
It starts growing a feasible alternating forest~$\calF$ rooted at the source nodes (line~\ref{line:forest_init}).
During the course of the execution, \textsc{Fast-Augment} \emph{scans} the arcs of $G$
in a certain order. For each node $u$ in~$G$ we keep the list $L(u)$ of all unscanned arcs leaving~$u$.
The following invariant is maintained:
\begin{numitem}
\label{eq:scanned_invariant}
    if $a = (u,v)$ is a scanned arc then either $a \in M$ or
    both $u^1$ and $v^2$ are $\calF$-reachable.
\end{numitem}

Consider an $\calF$-reachable node~$u^1$.
To enumerate the arcs leaving~$u^1$ in $\vec G(M)$ we fetch an unscanned arc $a = (u,v)$ from $L(u)$.
If $a \in M$ or $v^2$ is $\calF$-reachable then $a$ is skipped
and another arc is fetched.
(In the former case $a$ does not generate an arc leaving $u^1$ in $\vec G(M)$,
in the latter case $v^2$ is already $\calF$-reachable so $a$ can be made scanned
according to \refeq{scanned_invariant}.)

Otherwise, consider the arc $a_1 := (u^1,v^2)$ in $\vec G(M)$
and let $P_0$ denote the even alteranting feasible path from a root of $\calF$ to $u^1$.
Note that each node $x^2$ in $\vec G(M)$ (for $x \in VG$) is either a sink or has the unique arc leaving it.
A \emph{single} step occurs when $v^2$ is a sink (lines~\ref{line:single_step_first}--\ref{line:single_step_last}).
The algorithm constructs an augmenting path~$P_1 = P_0 \circ a_1$ leading to~$v^2$.
(Here $L_1 \circ L_2$ stands for the concatenation of $L_1$ and $L_2$.)
If $P_1$ is feasible, the current even factor gets augmented according to~\refclaim{aug_path} and \textsc{Fast-Augment} terminates.
Otherwise, forest growing stops and the algorithm proceeds to line~\ref{line:contraction_start} to deal with a contraction.

A \emph{double} step is executed when $v^2$ is not a sink (lines~\ref{line:double_step_first}--\ref{line:double_step_last}).
To keep the leafs of $\calF$ on even distances from roots, $\calF$ it must be extended by adding pairs of arcs.
Namely, there is a unique arc leaving $v^2$ in $\vec G(M)$, say $a_2 = (v^2,w^1)$ (evidently $(w,v) \in M$).
Moreover, $w^1$~is not a source node and $(v^2,w^1)$ is the only arc entering $w^1$.
Hence, $w^1$ is not $\calF$-reachable.
If $P_1 := P_0 \circ a_1 \circ a_2$ is feasible then $a_1$ and $a_2$
are added to $\calF$ thus making $v^2$ and $w^1$ $\calF$-reachable.
Otherwise, a contraction is necessary.

Now we explain how the algorithm deals with contractions at line~\ref{line:contraction_start}.
One has an even alternating feasible path~$P_0$ and an infeasible augmenting or even alternating
path~$P_1$ (obtained by extending $P_0$ by one or two arcs).
Put $M_0 := M \symdiff A(P_0)$ and $M_1 := M \symdiff A(P_1)$.
Let $C$ denote the unique odd cycle in $M_1$.
Put $G' := G / C$, $M' := M_0 / C$.
If $\textit{sparsify} = \false$ then \textsc{Fast-Augment} acts similar to \textsc{Simple-Augment},
namely, it makes a recursive call passing $G'$ and $M'$ as an input and, hence,
restarting the whole path computation.

Next, suppose $\textit{sparsify} = \true$.
In this case the algorithm tries to \emph{recover} some proper
alternating forest for the contracted digraph~$G'$ and the updated even factor~$M'$.
To accomplish this, a sparse digraph $H'$ is constructed and \textsc{Fast-Augment} is
recursively called for it (with $\textit{sparsify} = \false$).
The latter nested call may obtain a breakthrough, that is, find an even factor of smaller deficiency.
In this case, the outer call terminates immediately.
Otherwise, the nested call returns a complete proper alternating forest~$\bar \calF$ for an even factor $\bar M'$
in a digraph $\bar H'$ (obtained from $H'$ by contractions).
This forest is used by the outer call to continue the path-searching process.
It turns out that almost all of the arcs that were earlier fetched by the outer call
need no additional processing and may remain scanned w.r.t. the new, recovered forest.
This way, the algorithm amortizes arc scans during the outer call.

More formally, the algorithm constructs~$H'$ as follows.
First, take the node set of $G$, add all the arcs of~$M$ and all the arcs $(u, v) \in AG$ such that $(u^1, v^2)$ is present in~$\calF$.
Denote the resulting digraph by~$H$.
We need to ensure that $H$ is odd-cycle symmetric: if some arc $(u,v)$ in already added to~$H$
and the reverse arc $(v,u)$ exists in $G$ then add $(v,u)$ to~$H$.
Next, put $H' := H / C$.
Note that $H'$ is a sparse spanning subgraph of~$G'$ (i.e. $VH' = VG'$, $\abs{AH'} = O(n)$)
and $M'$ is an even factor in~$H'$.

The algorithm makes a recursive call $\textsc{Fast-Augment}(H',M',\false)$.
Let $\bar H'$ and $\bar M'$ be the resulting digraph and the even factor, respectively.
Compare the node sets of $\bar H'$ and $H$.
Clearly, $\bar H'$ is obtained from $H$ by a number of contractions ($G/C$ being one of them).
Let $Z_1, \ldots, Z_k$ be the maximum disjoint subsets of $G$ such that $\bar H' = G / Z_1 / \ldots / Z_k$.
Also, let $z_1, \ldots, z_k$ be the composite nodes in~$\bar H'$ corresponding to $Z_1, \ldots, Z_k$.
The algorithm applies these contractions to $G$ and constructs the digraph $\bar G' := G / Z_1 / \ldots / Z_k$.
Clearly $\bar M'$ is an even factor in both $\bar G'$ and $\bar H'$.
If $\mydef{\bar H',\bar M'} < \mydef{H',M'} = \mydef{G,M}$, then one has a breakthrough,
\textsc{Fast-Augment} terminates yielding $\bar G'$ and $\bar M'$.

Otherwise, the recursive call in line~\ref{line:recurse} also returns
a proper complete forest~$\bar \calF$ for $\bar H'$ and $\bar M'$.
Recall that some arcs in $G$ are marked as \emph{scanned}.
Since we identify the arcs of $\bar G'$ with their pre-images in $G$, one may speak of scanned arcs of~$\bar G'$.
The algorithm ``unscans'' certain arcs $a = (u, v) \in A\bar G'$
by adding them back to their corresponding lists $L(u)$ to ensure~\refeq{scanned_invariant}.
Namely, the arcs that belong to $M$ and are present in $\bar G'$ and
the arcs that enter any of the complex nodes $z_1, \ldots, z_k$ in $\bar G'$ are unscanned.
After this, the algorithm puts $G := \bar G'$, $M := \bar M'$, $\calF := \bar \calF$ and
proceeds with growing~$\calF$ (using the adjusted set of scanned arcs).

Finally, if \textsc{Fast-Augment} has scanned all the arcs of $G$ and
is unable to reach a sink, the resulting forest~$\calF$ is both proper and complete.
In particular, by \refeq{scanned_invariant} no augmenting path for $M$ exists.
By \refclaim{contraction_opt} this implies the maximality
of~$M$ in~$G$. The algorithm returns the current digraph~$G$, the current (maximum) even factor~$M$, and also
the forest~$\calF$, which certifies the maximality of~$M$.

The correctness of \textsc{Fast-Augment} is evident except for the case
when it tries to recover~$\calF$ and alters the set of scanned arcs.
One has to prove that \refeq{scanned_invariant} holds
for the updated forest and the updated set of the scanned arcs.
The proof of this statement is given in \refsec{correctness}.

\subsection{Complexity}

We employ arc lists to represent digraphs.
When a subset $U$ in a digraph~$\Gamma$ is contracted we enumerate all arcs incident to $U$ and
update the lists accordingly.
If a pair of parallel arcs appears after contraction, these arcs are merged,
so all our digraphs remain simple.
The above contraction of $U$ takes $O(\abs{V\Gamma} \cdot \abs{U})$ time.
During \textsc{Fast-Augment} the sum of sizes of the contracted subsets telescopes to $O(n)$,
so graph contractions take $O(n^2)$ time in total.
The usual bookkeeping allows to undo contractions and recover a maximum even factor
in the original digraph in $O(m)$ time.

Consider an invocation $\textsc{Fast-Augment}(\Gamma,N,\false)$ and let us bound its complexity
(including the recursive calls).
The outer loop of the algorithm (lines~\ref{line:loop_first}--\ref{line:loop_last}) enumerates the unscanned arcs.
Since $\textit{sparsify} = \false$, each arc can be scanned at most once,
so the bound of $O(\abs{A\Gamma})$ for the number of arc scans follows.
Using the appropriate data structures to represent even factors
the reachability checks in lines~\ref{line:feasibility_check_1} and~\ref{line:feasibility_check_2}
can be carried out in $O(\log{\abs{V\Gamma}})$ time (see~\refsec{feasibility_check} for more details).
Constructing $M_0$, $M_1$, and $C$ takes $O(\abs{V\Gamma})$ time.
This way, $\textsc{Fast-Augment}(\Gamma,N,\false)$ takes $O((k + 1) \abs{A\Gamma} \log \abs{V\Gamma})$ time,
where $k$ denotes the number of graph contractions performed during the invocation.

Next, we focus on $\textsc{Fast-Augment}(\Gamma,N,\true)$ call.
Now one may need to perform more than $\abs{A\Gamma}$ arc scans
since forest recovery may produce new unscanned arcs (line~\ref{line:unscan}).
Note that forest recovery totally occurs $O(\abs{V\Gamma})$ times (since each such occurrence leads to a contraction).
During each recovery $M$ generates $O(\abs{V\Gamma})$ unscanned arcs
or, in total, $O(\abs{V\Gamma}^2)$ such arcs for the duration of \textsc{Fast-Augment}.
Also, each node $z_i$ generates $O(\abs{V\Gamma})$ unscanned arcs
(recall that we merge parallel arcs and keep the current digraph simple).
The total number of these nodes processed during \textsc{Fast-Augment} is $O(\abs{V\Gamma})$
(since each such node corresponds to a contraction).
Totally these nodes produce $O(\abs{V\Gamma}^2)$ unscanned arcs.
Hence, the total number of arc scans is $O(\abs{A\Gamma} + \abs{V\Gamma}^2) = O(\abs{V\Gamma}^2)$.

Each feasibility check costs $O(\log \abs{V\Gamma})$ time, or $O(\abs{V\Gamma}^2 \log \abs{V\Gamma})$ in total.
Finally, we must account for the time spent in the recursive invocations during $\textsc{Fast-Augment}(\Gamma,N,\true)$.
Each such invocation deals with a sparse digraph and hence takes $O((k + 1) \abs{V\Gamma} \log \abs{V\Gamma})$ time
(where, as earlier, $k$ denotes the number of contractions performed by the recursive invocation).
Since the total number of contractions is $O(\abs{V\Gamma})$, the sum over all recursive invocations
telescopes to $O(\abs{V\Gamma}^2 \log \abs{V\Gamma})$.

The total time bound for $\textsc{Fast-Augment}(\Gamma,N,\true)$ (including the recursive calls)
is also $O(\abs{V\Gamma}^2 \log \abs{V\Gamma})$.
Therefore a maximum even factor in an odd-cycle symmetric digraph can be found in $O(n^3 \log n)$ time,
as claimed.


\nocite{*}
\bibliographystyle{plain}
\bibliography{main}

\begin{thebibliography}{10}

\bibitem{BP-89}
E.~Balas and W.~Pulleyblank.
\newblock The perfectly matchable subgraph polytope of an arbitrary graph.
\newblock {\em Combinatorica}, 9:321--337, 1989.

\bibitem{CSRL-01}
T.~Cormen, C.~Stein, R.~Rivest, and C.~Leiserson.
\newblock {\em Introduction to Algorithms}.
\newblock McGraw-Hill Higher Education, 2001.

\bibitem{cun-02}
W.~H. Cunningham.
\newblock Matching, matroids, and extensions.
\newblock {\em Mathematical Programming}, 91(3):515--542, 2002.

\bibitem{CG-97}
W.~H. Cunningham and J.~F. Geelen.
\newblock The optimal path-matching problem.
\newblock {\em Combinatorica}, 17:315--337, 1997.

\bibitem{CG-00}
W.~H. Cunningham and J.~F. Geelen.
\newblock Combinatorial algorithms for path-matching, 2000.
\newblock Manuscript.

\bibitem{CG-01}
W.~H. Cunningham and J.~F. Geelen.
\newblock Vertex-disjoint dipaths and even dicircuits, 2001.
\newblock Manuscript.

\bibitem{HK-73}
J.~E. Hopcroft and R.~M. Karp.
\newblock An {$n^{5/2}$} algorithm for maximum matchings in bipartite graphs.
\newblock {\em SIAM Journal on Computing}, 2(4):225--231, 1973.

\bibitem{IT-07}
S.~Iwata and K.~Takazawa.
\newblock The independent even factor problem.
\newblock In {\em Proceeinds of the 18th Annual ACM-SIAM Symposium on Discrete
  algorithms}, pages 1171--1180, 2007.

\bibitem{MV-80}
S.~Micali and V.~Vazirani.
\newblock An {$O(\sqrt{\abs{V}} \cdot \abs{E})$} algorithm for finding maximum
  matching in general graphs.
\newblock {\em Proc. 45st IEEE Symp. Foundations of Computer Science}, pages
  248--255, 1980.

\bibitem{pap-05}
G.~Pap.
\newblock A combinatorial algorithm to find a maximum even factor.
\newblock In {\em IPCO}, pages 66--80, 2005.

\bibitem{pap-07}
G.~Pap.
\newblock Combinatorial algorithms for matchings, even factors and square-free
  2-factors.
\newblock {\em Math. Program.}, 110(1):57--69, 2007.

\bibitem{PS-04}
G.~Pap and L.~Szeg\"{o}.
\newblock On the maximum even factor in weakly symmetric graphs.
\newblock {\em J. Comb. Theory Ser. B}, 91(2):201--213, 2004.

\bibitem{SW-02}
B.~Spille and R.~Weismantel.
\newblock A generalization of {E}dmonds' matching and matroid intersection
  algorithms.
\newblock In {\em Proceedings of the 9th International IPCO Conference on
  Integer Programming and Combinatorial Optimization}, pages 9--20, 2002.

\bibitem{T-08}
K.~Takazawa.
\newblock A weighted even factor algorithm.
\newblock {\em Mathematical Programming: Series A and B}, 115(2):223--237,
  2008.

\bibitem{tar-83}
R.~Tarjan.
\newblock {\em Data structures and network algorithms}.
\newblock Society for Industrial and Applied Mathematics, Philadelphia, PA,
  USA, 1983.

\bibitem{KT-09}
K.~Takazawa Y.~Kobayashi.
\newblock Even factors, jump systems, and discrete convexity.
\newblock {\em J. Comb. Theory Ser. B}, 99(1), 2009.

\end{thebibliography}

\newpage

\appendix

\noindent{\bf \Large Appendix}

\section{Correctness}
\label{sec:correctness}

In order to establish the correctness of the algorithm one needs to prove that
once the forest~$\calF$ gets recovered at line~\ref{line:recovered} of \textsc{Fast-Augment}
the set of scanned arc obeys property~\refeq{scanned_invariant}.
The latter is equivalent to the following:

\begin{lemma}
\label{lm:unscan_is_correct}
    Consider line~\ref{line:recovered} in \refalgo{fast_augment} and let $a = (u, v)$ be an arbitrary arc of $\bar G'$.
    Then either $a$ is unscanned or $a \in \bar M'$ or both $u^1$ and $v^2$ are $\bar \calF$-reachable.
\end{lemma}

First, we shall need a more convenient characterization of reachable nodes.
Consider a digraph~$\Gamma$ and a node $v \in V\Gamma$.
Construct a new odd-cycle symmetric digraph $\Gamma \ast v$ from $\Gamma$
by adding a new node $v'$ and an arc $(v,v')$.
\begin{lemma}
\label{lm:def_changes}
    For each $\Gamma$ and $v$ either $\mydef{\Gamma \ast v} = \mydef{\Gamma}$ or $\mydef{\Gamma \ast v} = \mydef{\Gamma} + 1$.
\end{lemma}
\begin{proof}
    Each even factor~$N$ in $\Gamma$ is also an even factor in $\Gamma \ast v$, hence
    $\mydef{\Gamma \ast v} \le \mydef{\Gamma} + 1$. Also, for an even factor $N^*$ in $\Gamma \ast v$
    put $N := N^* \setminus \set{(v,v')}$. Then, $N$ is an even factor in $\Gamma$ obeying
    $\abs{N} \ge \abs{N^*} - 1$. This implies $\mydef{\Gamma} \le \mydef{\Gamma \ast v}$, as required.
\end{proof}

\begin{lemma}
\label{lm:reachable}
    Let $N$ be a maximum even factor in $\Gamma$ and $\calT$ be a feasible alternating forest for $N$.
    If $v^1$ is $\calT$-reachable then $\mydef{\Gamma \ast v} = \mydef{\Gamma}$.
\end{lemma}
\begin{proof}
    Consider the even alternating path~$P$ from a root of $\calT$ to $v^1$.
    Since $\calT$ is feasible, $N' := N \symdiff A(P)$ is an even factor in~$\Gamma$ and no arc of $N'$ leaves~$v'$.
    Now $N^* := N' \cup \set{(v,v')}$ is an even factor in $\Gamma \ast v$.
    Therefore, $\mydef{\Gamma \ast v} \le \mydef{\Gamma \ast v, N^*} = \mydef{\Gamma,N} = \mydef{\Gamma}$
    implying $\mydef{\Gamma \ast v} = \mydef{\Gamma}$ by~\reflm{def_changes}.
\end{proof}

\begin{lemma}
\label{lm:not_reachable}
    Let $N$ be a maximum even factor in $\Gamma$ and $\calT$ be a complete alternating forest for $N$.
    If $v^1$ is not $\calT$-reachable then $\mydef{\Gamma \ast v} = \mydef{\Gamma} + 1$ and $N$ is a maximum
    even factor in $\Gamma \ast v$.
\end{lemma}
\begin{proof}
    The alternating forest~$\calT^*$ obtained from $\calT$ by adding a new root node~$(v')^1$ is complete for~$\Gamma \ast v$ and~$N$.
    Hence, $N$ is maximum in both $\Gamma$ and $\Gamma \ast v$ and $\mydef{\Gamma \ast v} = \mydef{\Gamma}$.
\end{proof}
Consider a node $v \in VG$. If $v \in VG - (Z_1 \cup \ldots \cup Z_k)$ then we
say that $v$ \emph{survives the contractions}. These nodes are both present in $G$ and $\bar G'$.

\begin{lemma}
\label{lm:v1_reachable}
    Suppose that a node $v \in VG$ survives the contractions and $v^1$ is $\calF$-reachable.
    Then $v^1$ is $\bar \calF$-reachable.
    Also, the nodes~$z_1^1, \ldots, z_k^1$ are $\bar \calF$-reachable.
\end{lemma}
\begin{proof}
    The transformation of~$H$ to~$\bar H'$ and of~$M$ to~$\bar M'$ may be viewed as follows:
    \begin{equation}
    \label{eq:contractions}
        \begin{array}{ccccccccccc}
            (H,M)               && (H,M_0)               && (H',M')                                                 &&&&&& (\bar H', \bar M') \\
            \rotatebox{90}{=}   && \rotatebox{90}{=}     && \rotatebox{90}{=}                                       &&&&&& \rotatebox{90}{=} \\
            (\Gamma^0,N^0) & \to & (\Gamma^0,N_0^0) & \to & (\Gamma^1,N^1) & \to & (\Gamma^1,N_0^1) & \to & \ldots & \to & (\Gamma^s,N^s)\\
        \end{array}
    \end{equation}
    Here $\Gamma^0, \ldots, \Gamma^s$ are odd-cycle symmetric digraphs,
    $N^i$, $N_0^i$ are even factors in~$\Gamma^i$ obeying $\abs{N^i} = \abs{N_0^i}$.
    Each $N_0^i$ fits some odd cycle $C^i$ in $\Gamma^i$ and $\Gamma^{i+1} = \Gamma^i / C^i$, $N^{i+1} = N_0^i / C^i$.

    Recall that the nested call to \textsc{Fast-Augment} in line~\ref{line:recurse} did not result into a breakthrough,
    so $\bar M'$ is maximum in $\bar H'$ and $\bar\calF$ is a complete alternating forest.

    Let us prove the first claim of the lemma.
    Consider a node~$v$ surviving the contractions such that $v^1$ is $\calF$-reachable in~$G$.
    Suppose towards a contradiction that $v^1$ is not $\bar\calF$-reachable in~$\bar H'$.
    By \reflm{not_reachable}, $\bar M' = N^s$ is a maximum even factor in $\bar H' \ast v = \Gamma^s \ast v$.
    Then, by \refclaim{contraction_opt},  $N_0^{s-1}$ is a maximum even factor in $\Gamma^{s-1} \ast v$
    and hence so is $N^{s-1}$ (by the equality of sizes). Proceeding this way in the backward direction
    we conclude that $N^i$ is a maximum even factor in $\Gamma^i \ast v$ for all $i = 0, \ldots s$.
    In particular, $N^0 = M$ is maximum in $\Gamma^0 = H$.
    Since $(v,v') \notin M$, $M$ is also maximum in~$H$ and $\mydef{H \ast v,M} = \mydef{H,M} + 1$.
    This contradicts \reflm{reachable} and the fact that $v^1$ is $\calF$-reachable.

    For the second claim, fix a node $v = z_i$ and suppose that $v^1$ is not $\bar \calF$-reachable in~$\bar H'$.
    Consider the sequence of transformations~\refeq{contractions} and
    suppose that $v$ was created as a complex node in $\Gamma^j$ while contracting an odd cycle $C^{j-1}$ in $\Gamma^{j-1}$.
    As indicated above, $N^j$ is a maximum even factor in $\Gamma^j \ast v$.
    The latter, however, is false since $N_0^{j-1}$ fits $C^{j-1}$ and, therefore, $N^j$ has no arcs leaving~$v$
    (cf. \reffig{contraction}(b) for an example).
    Hence, $N^j$ can be enlarged to $N^j \cup \set{(v,v')}$~--- a contradiction.
\end{proof}

\begin{lemma}
\label{lm:v2_reachable}
    Suppose that a node $v \in VG$ survives the contractions and $v^2$ is $\calF$-reachable.
    Then $v^2$ is $\bar \calF$-reachable.
\end{lemma}
\begin{proof}
    The node~$v^2$ cannot be a source, hence it is reached by some arc $(u^1,v^2)$ in~$\calF$,
    where $u^1$ is $\calF$-reachable and $a = (u,v) \notin M$.
    Let $a_0 = (u_0, v)$ be the image of~$a$ under contractions.
    Note that $a \in AH$ and $a_0 \in A\bar H'$.
    By \reflm{v1_reachable}, $u_0^1$ is $\bar \calF$-reachable.
    Therefore, if $a_0 \notin \bar M'$ then the completeness of $\bar \calF$ implies that $v^2$ is $\bar \calF$-reachable.
    It remains to consider the case $a_0 \in \bar M'$.
    The node $u_0^1$ is not a source (since $a_0 \in \bar M'$ leaves~$u_0$) but is $\bar \calF$-reachable.
    In the auxiliary bipartite digraph $u_0^1$ is entered by exactly one arc, namely $(v^2, u_0^1)$.
    Hence, $v^2$ must be $\bar \calF$-reachable, as required.
\end{proof}

Finally, we present the desired correctness proof.
\begin{proofof}{\reflm{unscan_is_correct}}
    Consider an arc $a = (u,v) \in A\bar G'$.
    If $a$ is not scanned then we are done.
    Otherwise let $a_0 = (u_0, v_0)$ denote its pre-image in~$G$.
    Here $u = u_0$ if $u_0$ survives the contractions or $u = z_i$ if $u_0 \in Z_i$.
    Also, since all arcs entering the complex nodes $z_1, \ldots, z_k$ are unscanned (line~\ref{line:unscan}),
    $v_0$ survives the contractions and hence $v = v_0$.
    Since the algorithm only decreases the set of scanned arcs, $a_0$ must also be scanned in~$G$.
    Clearly $a_0 \notin M$ since all the arcs that belong to~$M$ and are present in $\bar G'$ were as unscanned.
    Therefore, both $u_0$ and $v_0 = v$ are $\calF$-reachable in $G$ by the invariant~\refeq{scanned_invariant}.
    Applying \reflm{v1_reachable} to~$u_0$ and \reflm{v2_reachable} to~$v_0$
    we see that both $u$ and $v$ must be $\bar \calF$-reachable.
    The proof is now complete.
\end{proofof}

\section{Feasibility Checks}
\label{sec:feasibility_check}

This section explains how path feasibility checks, which are performed by \textsc{Fast-Augment}
at lines~\ref{line:feasibility_check_1} and~\ref{line:feasibility_check_2}, can be made efficient.
That is, given an even factor $N$ and a feasible even alternating path~$P_0$ we need to verify
that an even alternating or an augmenting path~$P_1$ (obtained from $P_0$ by appending one or two arcs) is feasible.
We make use of a certain data structure $\calD$ that maintains an even factor $M \symdiff A(P_0)$
as a collection of node-disjoint paths and cycles.
The following operations are supported by $\calD$:
\begin{itemize} \compact
    \item $\textsc{Insert}(u,v)$: assuming that $u$ is the end node of some path $P_u$ in $\calD$ and
    $v$ is the start node of some path $P_v$ in $\calD$, add the arc $(u,v)$ thus linking $P_u$ and $P_v$
    or, in case $P_u = P_v$, turning this path into a cycle;
    \item $\textsc{Remove}(u,v)$: assuming that $a = (u,v)$ is an arc belonging to some path or cycle in~$\calD$,
    remove $a$ thus splitting the path into two parts or turning the cycle into a path.
    \item $\textsc{Is-Odd-Cycle}(u,v)$: assuming that $a = (u,v)$ is an arc belonging to some path or cycle in~$\calD$,
    check if $a$ belongs to an odd cycle.
\end{itemize}

We make use of balanced search trees additionally augmented with \textsc{Split} and \textsc{Concatenate}
operations (e.g. red-black trees, see~\cite{CSRL-01,tar-83}) to represent paths and cycles in $\calD$.
This way, \textsc{Insert}, \textsc{Remove}, and \textsc{Is-Odd-Cycle} take $O(\log \abs{V\Gamma})$ time each.
Now checking if $P_1$ for feasibility is done by calling $\textsc{Insert}(u,v)$ and, in case of a double step,
$\textsc{Remove}(w,v)$, and finally making $\textsc{Is-Odd-Cycle}(u,v)$ request.
If the latter indicates that $P_1$ is not feasible, the changes in~$\calD$ are rolled back.

\medskip

During a $\textsc{Fast-Augment}(\Gamma,N,\false)$ call
we grow $\calF$ in a depth-first fashion and maintain the structure~$\calD$
corresponding to the current $\calF$-reachable node~$u^1$
(i.e. $\calD$ keeps the decomposition of $N \symdiff A(P)$
where $P$ is the path in $\calF$ from a root to~$u^1$).
When $\calF$ gets extended by arcs $(u^1,v^2)$ and $(v^2,w^1)$,
$w^1$ becomes the new current node and $\calD$ is updated accordingly
by the above $\textsc{Insert}(u,v)$ and $\textsc{Remove}(w,v)$ calls.
When the algorithm backtracks from $w^1$ to $u^1$, changes in $\calD$ are reverted.
This way, each feasibility check costs $O(\log \abs{V\Gamma})$ time

\medskip

Next, consider a $\textsc{Fast-Augment}(\Gamma,N,\true)$ call.
The above time bound of $O(\log \abs{V\Gamma})$ per check is only valid if we grow $\calF$ from scratch.
However, the algorithm also reuses the forest that is returned by the
nested \textsc{Fast-Augment} call in line~\ref{line:recurse}.
This incurs an overhead of $O(\abs{V\Gamma} \log \abs{V\Gamma})$ per forest recovery
(this additional time is needed to traverse the arcs that are present in the recovered forest~$\calF$ and update $\calD$ accordingly).
There are $O(\abs{V\Gamma})$ forest recoveries during the call
so the total overhead is $O(\abs{V\Gamma}^2 \log \abs{V\Gamma})$ time.
This does not affect the time bound of $\textsc{Fast-Augment}(\Gamma,N,\true)$.

\end{document}